
\input gtmacros
\input epsf

\input gtmonout
\volumenumber{2}
\volumeyear{1999}
\volumename{Proceedings of the Kirbyfest}
\pagenumbers{335}{342}
\papernumber{18}
\received{13 August 1998}\revised{26 February 1999}
\published{20 November 1999}

\def\zz{{\Bbb Z}}	
\def\qq{{\Bbb Q}}	
\def\calc{{\cal C}}
 
\def\cala{{\cal A}}

\reflist

\key{AK} {\bf S  Akbulut}, {\bf R Kirby}, {\it Branched
covers of surfaces in
$4$--manifolds}, Math. Ann. 252 (1979/80) 111--131.

\key{CG1} {\bf A Casson}, {\bf C Gordon}, {\it Cobordism
of classical knots,} from: ``A la recherche de la Topologie
perdue'', ed. by Guillou and Marin,  Progress in
Mathematics, Volume 62 (1986)  Originally published as
Orsay Preprint (1975)

\key{CG2} {\bf A Casson}, {\bf C Gordon}, {\it On slice
knots in dimension three}, Proc. Symp. Pure Math. 32
(1978) 39--54

\key{D} {\bf S Donaldson}, {\it An application of gauge
theory to four--dimensional topology}, J. Differential Geom.
18 (1983) 279--315

\key{F} {\bf R Fox} {\it A quick trip through knot theory}, 
from: ``Topology of 3--manifolds and related topics'', (Proc. The
Univ. of Georgia Institute, 1961)
Prentice--Hall, Englewood Cliffs, NJ (1962) 120--167 

\key{FM} {\bf R Fox}, {\bf J Milnor}, {\it Singularities of
$2$--spheres in
$4$--space and cobordism of knots}, Osaka J. Math. 3 (1966) 257--267 

\key{Fr} {\bf M Freedman}, {\it The topology of
four--dimensional manifolds}, J. Diff. Geom. 17
(1982) 357--453 

\key{FQ} {\bf M Freedman}, {\bf F Quinn}, {\it Topology of
4--manifolds}, Princeton Mathematical Series, 39, Princeton
University Press, Princeton, NJ (1990)

\key{Gi} {\bf P Gilmer},  {\it Slice knots in $S\sp{3}$}, 
Quart. J. Math. Oxford Ser. (2) 34 (1983)
305--322 

\key{Go} {\bf C Gordon}, {\it Problems}, from: ``Knot Theory'', ed.
J-C Hausmann,  Springer Lecture Notes no. 685 (1977)

\key{J} {\bf B Jiang}, {\it A simple proof that the
concordance group of algebraically slice knots is infinitely
generated}, Proc. Amer. Math. Soc. 83 (1981) 189--192

\key{K1} {\bf R Kirby}, {\it Problems in low dimensional
manifold theory}, in Algebraic and Geometric Topology
(Stanford, 1976), vol 32, part II of Proc. Sympos. Pure
Math. 273--312

\key{K2} {\bf R Kirby}, {\it Problems in low dimensional
manifold theory}, from: ``Geometric Topology'', AMS/IP Studies in
Advanced Mathematics, ed. W Kazez (1997)

\key {Le1} {\bf J Levine}, {\it Knot cobordism groups in
codimension two},  Comment. Math. Helv. 44 (1969) 229--244

 \key {Le2} {\bf J Levine}, {\it Invariants of knot
cobordism}, Invent. Math. 8 (1969) 98--110

\key{Li} {\bf R Litherland}, {\it Cobordism of Satellite
Knots}, from: ``Four--Manifold Theory'', Contemporary Mathematics,
eds. C Gordon and R Kirby, American Mathematical Society,
Providence RI (1984) 327--362

\key{LN} {\bf C Livingston}, {\bf S Naik}, {\it Obstructing
4--torsion in the classical knot concordance group},
to appear in Jour. Diff. Geom.

\key{M} {\bf K Murasugi}, {\it On a certain numerical
invariant of link types}, Trans. Amer. Math. Soc. 117
(1965) 387--422

\key{T} {\bf A Tristram}, {\it Some cobordism invariants for
links}, Proc. Camb. Phil. Soc. 66 (1969) 251--264
  
\endreflist

\title{Order 2 Algebraically Slice Knots}                    

\authors{Charles Livingston}                  
 
\address{Department of Mathematics, Indiana
University\\Bloomington, Indiana  47405, USA}

\email{livingst@indiana.edu}                     

\abstract  
The concordance group of algebraically slice knots is the subgroup of
the classical knot concordance group formed by algebraically slice
knots.  Results of Casson and Gordon and of Jiang showed that this group
contains in infinitely generated free (abelian) subgroup.  Here it is
shown that the concordance group of algebraically slice knots also
contain elements of finite order; in fact it contains an infinite
subgroup generated by elements of order 2.
\endabstract


\makeshorttitle

\section{Introduction}

The classical knot concordance group, $\calc$, was defined by
Fox [\ref{F}] in 1962.  The work of Fox and Milnor
[\ref{FM}], along with that of Murasugi [\ref{M}] and
Levine [\ref{Le1}, \ref{Le2}], revealed fundamental aspects
of the structure of
$\calc$. Since then there has been tremendous
progress in 3-- and 4--dimensional geometric topology, yet
nothing more is now known about the underlying group structure of $\calc$
than was known in 1969.  In this paper we will describe new and unexpected
classes of order 2 in $\calc$.

What is known about $\calc$ is quickly summarized.  It is a countable
abelian group.  According to [\ref{Le2}] there is a
surjective homomorphism of
$\phi \co \calc \rightarrow 
 \zz^\infty \oplus\zz_2^\infty \oplus\zz_4^\infty$.  The results of
[\ref{FM}] quickly yield an infinite set of elements of order 2
in
$\calc$, all of which are mapped to elements of order 2 by
$\phi$.  
	
The results just stated, and their algebraic consequences, present all
that is known concerning the purely algebraic structure of $\calc$ in
either the smooth or topological locally flat category.  For instance, one
can conclude that elements of order 2 detected by homomorphisms to $\zz_2$,
such as the Fox--Milnor examples, are not evenly divisible, but it remains
possible that any given countable abelian group is a subgroup of
$\calc$, including such groups as the infinite direct sum of copies of
$\qq$ and
$\qq /  \zz$.  Most succinctly, we know that $\calc$ is isomorphic to a
direct sum  $\zz^\infty \oplus\zz_2^\infty \oplus G$, but all that is
known about $G$ is that it is  countable and  abelian. In  particular, Fox
and Milnor's original question on the existence of torsion of order other
than 2 remains completely open.  Other basic questions regarding the
structure of
$\calc$ appear in [\ref{Go}, \ref{K1}, \ref{K2}].

More is known about the pair $(\calc, \phi)$.  Casson and
Gordon [\ref{CG1}, \ref{CG2}] showed that the kernel of
Levine's homomorphism, $\cala$ (the {\sl concordance group of
algebraically slice knots}), is nontrivial; Jiang showed that
Casson and Gordon's examples provide an infinitely generated
free subgroup of  $\cala$.  

(We should observe here that
Alexander polynomial 1 knots are known to represent classes in $\cala$. 
Freedman's work [\ref{Fr}] implies that all such knots are
topologically locally flat slice. Donaldson's work [\ref{D}]
implies that some such knots are not smoothly slice.  Needless
to say, the accomplishments of both [\ref{Fr}] and [\ref{D}]
have been revolutionary in the study of 4--manifolds.  However,
it is perhaps surprising that neither has revealed any further
group theoretic structure of either $\calc$ or the pair
$(\calc,
\phi)$.)

It is becoming clear that any unexpected complexity in $\calc$ will appear
in $\cala$, that is, among algebraically slice knots:  if any odd torsion
exist, it obviously must  be in $\cala$, and recent work
[\ref{LN}] showing that infinite collections of knots that map
to elements of order 4 under Levine's homomophism are of
infinite order in $\calc$ supports the conjecture that any
4--torsion must also be in $\cala$.

In this paper we construct an infinite family of order 2 elements in 
$\cala$.  These are the first such examples, and the first examples of
any type showing that $\cala$ has any structure beyond that demonstrated
by Jiang. Our methods  apply in the smooth setting, as did the
original work of Casson and Gordon, but  work  of Freedman
[\ref{FQ}] shows that they apply in the topological locally
flat category also. 

Thanks are due to Larry Taylor for pointing out the particular problem
being addressed here, and to Zhenghan Wang for observing a
simplification of our original construction. 

\section{Basic building blocks}

The basic idea of the construction of algebraically slice order 2
knots is to take the connected sums of pairs of  algebraically
concordant negative amphicheiral knots.  If the knots aren't concordant
then the connected sum will be of order 2 in $\cala$.  The trick
is to find an infinite collection of such pairs so that it is possible to
prove that they have the desired properties.  Our examples,
$ J_n$, will be built as connected sums of pairs of knots
$K_{T}$ that we first examine.  For an arbitrary knot,
$T$, consider the knot
$K_{T}$ illustrated in Figure 1 along with a surgery diagram of
$M_{T}$, the 2--fold branched cover of $S^3$ branched over $K_{T}$, drawn using
the algorithm of Akbulut and Kirby [\ref{AK}].    The
illustration indicates that
$K_{T}$ bounds a genus one Seifert surface so that one band in the surface
has the knot
$T$ tied in it and the other band has the mirror image of $T$, $-T$, tied
in it.  (More precisely, in Figure 1 the tangle $-T$ is obtained from the
tangle $T$ by changing all the crossings.)  The bands are twisted so that
$K_{T}$ has Seifert form 
$$\pmatrix{1&1\cr 0&-1}.$$

\figure
\epsfxsize=3.5truein
\epsfbox{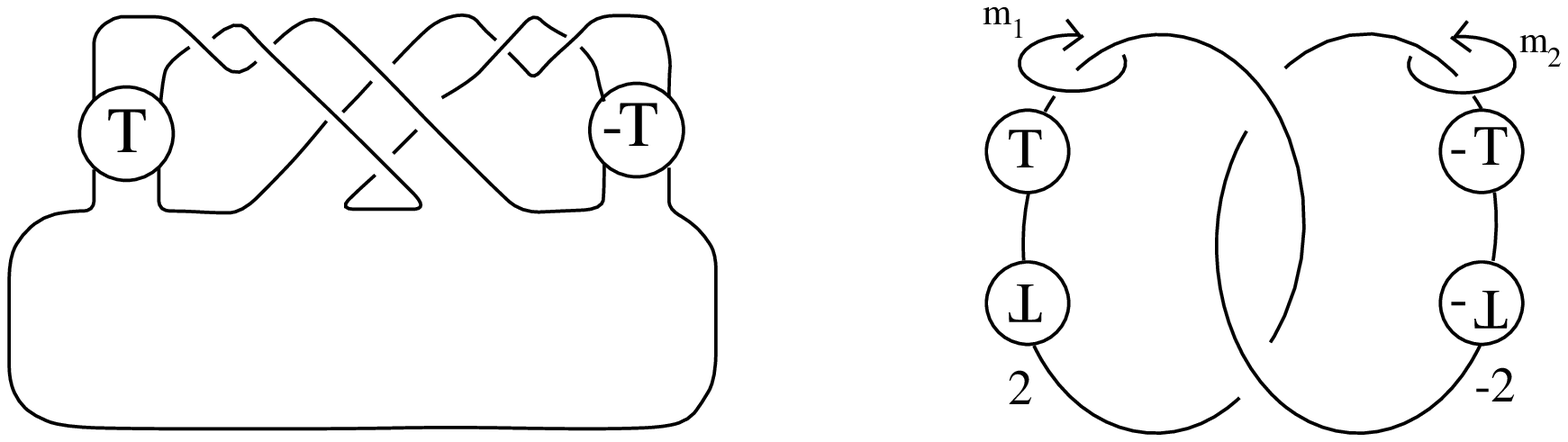}
\endfigure

\proc{Lemma}  The knot  $K_{T}$ is of order 2 in the knot
concordance group, $\calc$.\endproc 

\prf  Changing all the crossings in $K_{T}$ is easily seen
to have the effect of simply reversing its orientation; that is, $K_{T}$
is negative amphicheiral and is hence of order 1 or 2 in $\calc$.  (Stated
differently, $K_{T} = -K_{T}$, so $K_{T}$ is of order at most 2.)

The Alexander polynomial of $K_{T}$ is $t^2 - 3t + 1$, which is
irreducible, since the discriminant, $5$, is not a perfect square.
According to [\ref{FM}] this obstructs a knot from being
slice.  It follows that the order of $K_{T}$ is exactly 2. \endprf

We will need to understand the 2--fold branched cover of $S^3$, $M_{T}$.

\proc{Lemma}  $H_1(M_{T}) = \zz_{5}$ and is generated by the
meridian labeled $m_1$ in Figure 1.  Any homomorphism
$\phi \co H_1(M_{T}) \rightarrow \zz_{5}$ taking value $a$ on $m_1$ takes
value  $3a$ on $m_2$.\endproc

\prf  A relation matrix for the homology 
$H_1(M_{T})$ with respect to $m_1$ and $m_2$, computed using its surgery
presentation, is given by the  matrix  
$$\pmatrix{2&1\cr1&-2}.$$ This matrix presents the cyclic group $\zz_5$. 
The relations imply that $m_2$ = $-2m_1$, or $m_2 = 3m_1$.  \endprf

Our final calculation regarding $K_T$ is of its Casson--Gordon invariants. 
Associated to any representation $\chi \co H_1(M_{T}) \rightarrow \zz_5$
there is a rational invariant, denoted $\sigma_1(\tau (K_T,\chi))$ in
[\ref{CG1}].  We do not need the exact value of this invariant,
but only its dependence on $T$.  Let $K_0$ be the Figure 8
knot, obtained when $T$ is trivial.  Work of Litherland
[\ref{Li}] on computing $\tau$ for satellite knots yields the
desired result, Lemma 2.3 below.  Here, $\sigma_p$, $p \in \qq /
\zz$,  denotes the Tristram--Levine signature of a knot
[\ref{T}], given by the signature of the hermetian matrix
$$(1-e^{p2\pi i})V + (1-e^{-p2\pi i})V^t,$$ where $V$ denotes the
Seifert matrix of the knot.

\proc{Lemma}  If $\chi_a\co   H_1(M_{T}) \rightarrow \zz_5$ takes
value $a$ on $m_1$ then $$\sigma_1(\tau (K_T,\chi_a)) = \sigma_1(\tau
(K_0,\chi)) + 2\sigma_{a/5}(T) + 2\sigma_{3a/5}(-T).$$ \endproc

\prf The knot $K_T$ is a  satellite knot with
companion $T$, winding number 0, and with orbit $K'_{-T}$, the
knot formed from
$K_T$ by removing the knot $T$ in the left band.  Applying $\sigma_1$ to
the equation given as Corollary 2 in [\ref{Li}] gives:
$$\sigma_1(\tau (K_T,\chi_a)) =  \sigma_1(\tau
(K'_{-T},\chi_a)) + 2\sigma_{a/5}(T).$$  

Repeating this companionship argument to remove the $-T$ from the band in
$K'_{-T}$ yields the desired result.  The second signature is evaluated at
$3a/5$ because, according to Lemma 2.2, $\chi$ takes that value on $m_2$. 
(In Litherland's notation, $\chi(x_i) = 3a$, viewed as an element of
$\zz_5$.)
\endprf

\proclaim{Remark}\rm  We should remind the reader here that both the
Casson--Gordon invariant and the Tristram--Levine signature functions are
symmetric under a sign change; that is, $\sigma_1(\tau (K_T,\chi_a))
=\sigma_1(\tau (K_T,\chi_{-a}))$ and $\sigma_{p}(T) =\sigma_{-p}(T)$.
\endproc
 
Our ultimate examples, $J_i$, will be of the form $K_0 \# K_{T_i}$ for
particular $T_i$ which yield nontrivial signature values in the formula
of Lemma 2.3.  Explicit examples will be obtained by letting $T_i$  be 
connected sums of $(2,7)$--torus knots, so we conclude this section with the
following computation.

\proc{Lemma}  For the $T$ the $(2,7)$--torus knot and $a \ne 0 $
mod 5,    $\sigma_{a/5}(T) + \sigma_{3a/5}(-T) = 4 \hbox{ or -}4$ depending
on whether $a = \pm 2$ mod 5 or $a = \pm 1$ mod 5.\endproc
 
\prf  The signature function of $T$,
$\sigma_p$, is given (by definition) as the
signature of the form
$(1-e^{p2\pi i})V + (1-e^{-p2\pi i})V^t$ where  $V$ is a Seifert matrix for
$T$:
$$\pmatrix{1&1&0&0&0&0\cr 0&1&1&0&0&0\cr 0&0&1&1&0&0\cr 
0&0&0&1&1&0\cr 0&0&0&0&1&1\cr 0&0&0&0&0&1\cr}.$$ 
The signature of this
form is easily computed to be:

$$\sigma_{p}(T) = \cases{0, & if $0 < p <1/14$)
\cr 
         2, &if $1/14 < p < 3/14$ \cr
         4, &if $3/14 < p <  5/14$ \cr
         6, &if $5/14 < p < 7/14.$ \cr
}$$

The result now follows, since
one of the  signatures of $T$ that appears will be either 2 or 6 while the
signature of $-T$ that appears will be either $-6$ or $-2$, respectively.
\endprf

\section{Infinite 2--torsion among
algebraically slice knots}

\proclaim{Definition} Let $T$ represent the $(2,7)$--torus knot, let
$T_i = \#_iT$, let $K_i = K_{T_i}$ and let $J_i = K_{0} \#
K_{i}$.\endproc

If $i =0$ then  $\#_i T$ denotes the
unknot.  It follows that $J_0$ is the connected sum of the Figure 8 knot with
itself and is hence slice.  Also note that the definition makes sense for $i <
0$, letting $\#_i T$ denote the connected sum of $-i$ copies of the
$(2,-7)$--torus knots in that case.

\proc{Lemma} Each knot, $J_i$, is algebraically slice and of
order at most two in $\calc$ (and hence also in $\cala$). \endproc

\prf  Since each $K_{i}$ is of order two, the
connected sum of two of them is of order 1 or 2.  Also, since all the
$K_{i}$ have the same Seifert form, the Seifert form of an order two
knot, the connected sum of two of them is algebraically slice.\break\endprf

\proc{Theorem} For $i \ne j$, the knots $J_i$ and $J_j$ are not
concordant.\endproc 

\prf For arbitrary $i$ and $j$, if $J_i$ and $J_j$ are
concordant then $J_i \#J_j$ is slice.  Expanding, we have that $K_0 \# K_i
\# K_0 \#K_j$ is slice. Since $K_0$ is of order 2, this implies that $K_i
\# K_j$ is slice. 

The 2--fold branched cover of $S^3$ branched over $K_i \# K_j$ has first
homology splitting naturally as $\zz_5 \oplus \zz_5$.  There is a $\qq/\zz$
valued nonsingular symmetric linking form on this homology group.  Since the
covering space is a connected sum, the linking form splits along the given
direct sum decomposition.  Furthermore, the linking form takes the same value
on the merdians $m_1$ in each of the two individual branched covering spaces,
since the surgery matrices determine the linking form and are the same for
both covering spaces.

According to Casson and Gordon, if $K_i \# K_j$ is slice, there is some
vector $v \in \zz_5 \oplus \zz_5$ with self linking 0 so that for any $\zz_5$
valued character that vanishes on $v$ the associated Casson--Gordon
invariant must vanish.  With the given linking form, $v$ must be a
multiple of either $(2,1)$ or $(2,-1)$, and hence we can consider the
$\chi$ that takes value $1$ on the $m_1$ generator of $H_1(M_{T_i})$ and
value
 $\pm 2$ on the $m_1$ generator of $H_1(M_{T_j})$. 

Applying Lemma 2.3 along with the additivity of Casson--Gordon invariants
[\ref{Gi}] shows that if $K_i \# K_j$ is slice then
$$\eqalign{\sigma_1(\tau(M_{T_0},\chi_1)) + 2\sigma_{1/5}(T_i) +
2\sigma_{2/5}(-T_i) & +
\sigma_1(\tau(M_{T_0},\chi_2))\cr
& + 2\sigma_{2/5}(T_j)  + 2\sigma_{1/5}(-T_j)
=0.\cr}$$

Consider first the case that $i = 0 = j$.  In this case $K_i \# K_j$ is
slice, so the above formula yields that $$\sigma_1(\tau(M_{T_0},\chi_1)) +
\sigma_1(\tau(M_{T_0},\chi_2))
 =0.$$ Hence we can simplify the equation to find that if $K_i \# K_j$ is
slice, then $$ 2\sigma_{1/5}(T_i) + 2\sigma_{2/5}(-T_i) + 
2\sigma_{2/5}(T_j) + 2\sigma_{1/5}(-T_j) =0.$$ From Lemma 2.4 this
simplifies to be $8(j-i) =0$.  Clearly, if $i \ne j$ this yields a
contradiction. \endprf

\proc{Corollary} $\cala$ contains a subgroup isomorphic to
$\zz_2^\infty$. \qed \endproc\eject

\references
\Addresses\recd

\bye